\documentclass[preprint, 11pt]{elsarticle}
\usepackage[b5paper, margin={0.5in,0.65in}]{geometry}
\usepackage[T1]{fontenc}
\usepackage[latin9]{inputenc}
\usepackage{amsmath}
\usepackage{amsthm}
\usepackage{amssymb}
\usepackage{enumerate}   
\usepackage{color}
\usepackage{mathcommand}
\usepackage[%
  pagebackref, 
  breaklinks, 
  colorlinks=true, linkcolor=blue, citecolor=blue, urlcolor=blue
]{hyperref}

\numberwithin{equation}{section}

\def \1{{\bf 1}}
\def\C{{\mathbb{C}}}
\def\Z{{\mathbb{Z}}}

\def\N{{\mathbb{N}}}

\def\Y{{\mathcal{Y}}}
\def\I{{\mathcal{I}}}

\def\wt{{\rm wt}}
\def\w{{\omega}}

\def\i{{\mathrm{i}}} 
\def\dim{{\rm dim}}
\def\ker{{\rm ker}}

\def\Lz{{\text{Log}(z)}}

\def\id{{\rm id}}

\theoremstyle{definition}
\newtheorem{lemma}{Lemma}[section]
\newtheorem{theorem}[lemma]{Theorem}

\newtheorem{proposition}[lemma]{Proposition}
\newtheorem{definition}[lemma]{Definition}

\newtheorem{remark}[lemma]{Remark}
\newtheorem{mainthm}{Theorem}

\allowdisplaybreaks

\begin{document}
\begin{frontmatter}

\title{Cofiniteness for Twisted Fusion Products in Vertex Operator Algebra Theory}

\author[1]{{Chao} {Yang}}
\ead{yangcabc@163.com}

\author[2]{{Yiyi} {Zhu}}
\ead{chengdianyifan@126.com}

\affiliation[1]{
  organization={School of Mathematics, Southwest Jiaotong University}, 
  addressline={}, 
  city={Chengdu}, 
  postcode={611756}, 
  country={China},
}

\affiliation[2]{
  organization={School of Mathematics and Statistics, Guangdong University of Technology}, 
  addressline={}, 
  city={Guangzhou}, 
  postcode={510520}, 
  country={China},
}

\begin{abstract}
Let $V$ be a vertex operator algebra equipped with two commuting finite-order automorphisms $g_1$ and $g_2$, and set $g_3 = g_1 g_2$. For $k = 1, 2, 3$, let $W^k$ be a $g_k$-twisted $V$-module. Assuming that $W^1$ and $W^2$ are $C_1$-cofinite and that there exists a surjective twisted logarithmic intertwining operator of type $\binom{W^3}{W^1 \  W^2}$, 
we prove that $W^3$ is also $C_1$-cofinite. The cofiniteness follows from the finite-dimensionality of the solution space of an associated complex-coefficient linear differential equation. 
As an application,  under the condition of $C_1$-cofiniteness, we establish the finiteness of the fusion rules and construct the fusion product.
\end{abstract}


\end{frontmatter}

\tableofcontents

\section{Introduction}
The fusion product (or tensor product) theory provides a powerful framework for studying the representation theory of vertex operator algebras (VOAs). 
Compared to the classical case, the fusion product for VOAs is defined in terms of intertwining operators \cite{FHL93},
making its construction and associativity proof far from straightforward.
An important result in this theory states that, under suitable conditions, the representation category of a VOA forms a modular tensor category \cite{H08},
which has inspired a wealth of research at the intersection of VOAs and tensor category theory.

Given the fundamental role of twisted modules in vertex operator algebras, 
particularly in orbifold theory \cite{CM16,DLM98, DLM00, DNR25, DRX17},
a natural and important question is how to extend the theory of fusion products from untwisted to twisted modules. 
Beyond this extension itself, a key motivation is the prospect that such twisted fusion products would endow the category of twisted modules with a $G$-crossed tensor structure \cite{DH25}, thereby enabling a tensor-categorical approach to the representation theory of VOAs.
The main goal of this paper is to construct the fusion product of $C_1$-cofinite generalized twisted modules 
and to prove that this fusion product preserves $C_1$-cofiniteness.

The theory of twisted fusion products has been studied from various perspectives.
Xu first generalized intertwining operators to twisted modules in \cite{X95}, establishing a foundation for later developments. 
The concept of twisted fusion products was introduced by Dong, Li, and Mason in \cite{DLM96}, 
and they provided explicit constructions for certain twisted modules.
A more detailed study of twisted intertwining operators was given in \cite{LS23}. For strongly rational VOAs, the second author provided a construction of the fusion product for twisted modules in \cite{Z25}, extending Dong and Jiang's bimodule theory \cite{DJ08a, DJ08b}. Twisted fusion products have also been extensively studied in the context of permutation orbifolds in \cite{DLXY24, DXY25}. 
The twisted Verlinde formula was investigated in \cite{DLin25}.
These developments typically assume that the relevant automorphisms are commuting and of finite order. On the other hand, Huang introduced a more general notion of twisted intertwining operators for non-commuting automorphisms in \cite{H18}. Recently, \cite{DH25} constructed twisted fusion products under suitable assumptions.

This paper aims to extend Miyamoto's method \cite{M14} for constructing untwisted fusion products to the twisted case. This extension is particularly challenging because the Jacobi identity for twisted intertwining operators is much more complicated, making the derivation of key differential equations technically difficult (see Section~\ref{sec5}).
Our approach also differs from \cite{M14} in several technical aspects (see Section \ref{sec6-1}), allowing us to present a simpler and more direct construction of the twisted fusion product, particularly in the $C_2$-cofinite case (see Section \ref{sec6-2} ).

We now present some details. 
Let $V$ be a vertex operator algebra equipped with two commuting finite-order automorphisms $g_1$ and $g_2$, and set $g_3 = g_1 g_2$. 
Suppose the order of $g_k$ is $T_k$ for $k=1, 2, 3$. Let $W^1$ and $W^2$ be $C_1$-cofinite $g_1$- and $g_2$-twisted generalized $V$-modules, respectively. 
Let $\mathcal{F}_{(W^1, W^2)}$ be the set of pairs $(U, \mathcal{Y})$, where $U$ is a generalized $g_3$-twisted $V$-module and $\mathcal{Y} \in \mathcal{I} \binom{U}{W^1 \ W^2}$ is a surjective intertwining operator, modulo the equivalence relation defined by $(U, \mathcal{Y}_1) \cong (W, \mathcal{Y}_2)$ if there exists a $g_3$-twisted $V$-module isomorphism $f: U \to W$ satisfying $f \circ \mathcal{Y}_1 = \mathcal{Y}_2$.
Our first main result is:
\begin{mainthm}
For any $(U, \Y) \in \mathcal{F}_{(W^1,W^2)}$, we have
   \[
   \text{dim} \left(U/C_1(U)\right) \leq \text{dim} \left(W^1/C_1(W^1)\right) \text{dim} \left(W^2/C_1(W^2)\right).
   \]
   In particular, $U$ is $C_1$-cofinite.
\end{mainthm}

To construct the fusion product,  we introduce the following notation. Let 
\[
\mathcal{S} = \prod_{(U, \mathcal{Y}) \in \mathcal{F}(W^1,W^2)} U.
\]
Then, there is a natural action of $V$ on the product space $\mathcal{S}$ given by 
\[
Y_{\mathcal{S}}(-, x)=\prod_{(U, \mathcal{Y}) \in \mathcal{F}(W^1,W^2)}Y_{U}(-, x).
\]
For $w^1 \in W^1$, $w^2 \in W^2$, $n \in \mathbb{C}$, and $k \in \mathbb{N}$, we define
\[
(w^1)^{\diamond}_{n;k}w^2 = \prod_{(U, \mathcal{Y}) \in \mathcal{F}(W^1,W^2)} \big( (w^1)^{\mathcal{Y}}_{n;k}w^2 \big) \in \mathcal{S}.
\]
Let 
\[
W^1 \Diamond W^2 = \operatorname{span} \left\{ (w^1)^{\diamond}_{n;k}w^2 \ \middle| \ w^1 \in W^1, w^2 \in W^2, n \in \mathbb{C}, k \in \mathbb{N} \right\} \subset \mathcal{S}.
\]
Finally, define the operator
\[
\mathcal{Y}^{\diamond}(w^1,x)w^2 = \sum_{\substack{n \in \mathbb{C} \\ k \in \mathbb{N}}}
(w^1)^{\diamond}_{n;k}w^2 \, x^{-n-1} (\log x)^k \in (W^1 \Diamond W^2)[[\log x]]\{x\}.
\]
Then, our second main result is:
\begin{mainthm} \label{th-B}
 $ W^1 \Diamond W^2 \in \mathcal{F}_{(W^1,W^2)} $ and $(W^1 \Diamond W^2, \mathcal{Y}^{\diamond}) $  is the fusion product of $W^1$ and $W^2$.
\end{mainthm}

We note that if the tensor product of $W^1$ and $W^2$ exists, it is straightforward to verify that $W^1 \Diamond W^2$ is a generalized $g_3$-twisted $V$-module and that $\mathcal{Y}^{\diamond}$ defines a twisted intertwining operator. Moreover, $(W^1 \Diamond W^2, \mathcal{Y}^{\diamond})$ is isomorphic to the tensor product $(W^1 \boxtimes W^2, \mathcal{Y}^{\boxtimes})$.
This observation is precisely what motivates our construction of the fusion product in this manner.

However, without assuming the existence of fusion products, Theorem \ref{th-B} becomes non-trivial. 
Our proof relies essentially on the existence of complex numbers $\lambda_1, \dots, \lambda_n$ and a positive integer $d \in \mathbb{N}$ such that 
the following holds for any $(U, \Y) \in \mathcal{F}{(W^1,W^2)}$:
\begin{enumerate}[{(1)}]
\item \label{pro1}
The weight set $\wt(U)$ satisfies $\wt(U) \subset \bigcup_{i=1}^n \left( \lambda_i + \frac{1}{T_3} \mathbb{N} \right)$.
\item $\left(L(0)_n\right)^d U = 0$,
\end{enumerate}
where $L(0)_n$ is the nilpotent part of $L(0)$.
The crucial difference between our construction of the fusion product and that in \cite{M14} lies in our proof of the existence of a uniform bound $d$. 
This existence first guarantees that every $(w^1)^{\diamond}_{n;k}w^2 \in \mathcal{S}  $ is a generalized eigenvector of $L(0)$. 
Furthermore, it ensures that there exists a positive integer $K$ such that $\Y(w^1, x)w^2 \in W^3\{x\}[\text{log}x]_{\leq K} $ for 
any $(U, \Y) \in \mathcal{F}{(W^1,W^2)}$, which consequently implies that $\mathcal{Y}^{\diamond}(w^1,x)w^2 \in W^3\{x\}[\text{log}x]_{\leq K} $.
On the other hand, similar to \cite{M14}, the property (\ref{pro1}) above ensures that both $Y_{ \mathcal{S} }(v, x)s$ and $\mathcal{Y}^{\diamond}(w^1,x)w^2 $ are lower-truncated.

Finally, when $V$ is $C_2$-cofinite, we can straightforwardly construct fusion products.
Let $S^1, S^2, \cdots, S^p$ be a complete set of representatives of the isomorphism classes of irreducible grading restricted generalized $g_3$-twisted $V$-modules, with conformal weights $u_1, u_2, \cdots, u_p$, respectively.
For a grading restricted generalized $g_3$-twisted $V$-module $W$, define the natural number $d_W$ by 
\[
d_W=\text{dim} W_{(u_1)} + \cdots +\text{dim} W_{(u_p)}. 
\]
We show that $d_W$ is uniformly bounded for all $W\in (W, \Y)\in \mathcal{F}(W^1, W^2)$. Let $(W^1 \odot W^2, \Y^\odot) \in \mathcal{F}_{(W^1,W^2)}$ be such that $d_{W^1 \odot W^2}$ is maximal. Our third main result is:
\begin{mainthm}
   The pair $(W^1 \odot W^2, \Y^\odot)$ is the fusion product of $W^1$ and $W^2$.
\end{mainthm}

This paper is organized as follows: Section \ref{sec2} presents some basic results from analysis, which will be used later; In Section \ref{sec3}, we set up the stage by recalling and introducing some definitions, such as twisted module, $C_1$-cofiniteness, twisted intertwining operator, and twisted fusion product; In Section \ref{sec4}, we present our first main result, i.e., the $C_1$-cofiniteness of the fusion product of two $C_1$-cofinite generalized twisted modules; In Section \ref{sec5}, we provide constructions of fusion product of twisted modules to guarantee their existence. 

\section{Basic Facts from Analysis}\label{sec2}

Let
$$
\mathbb{C}_f\{x\} = \{ \sum_{n \in \C} a_n x^n \ |\ a_n \in \mathbb{C},\ a_n = 0 \text{ for all but finitely many } n \}.
$$

Let  
\[
\Omega = \mathbb{C} \setminus [0, \infty)
\]  
denote the complex plane slit along the positive real axis. Then $\Omega$ is a simply connected domain in $\mathbb{C} \setminus \{0\}$.  
In this paper, we consistently fix the branch of the logarithm on $\Omega$ as follows:  
\[
\operatorname{Log} z = \ln |z| + i \operatorname{Arg} z, \quad \text{where } 0 < \operatorname{Arg} z < 2\pi.  
\]

For $\alpha \in  \C  $, we define the power function  by
\[
z^{\alpha} = e^{\alpha \Lz}, \quad \text{for } z \in \Omega.
\]
It is worth noting that both $\Lz$ and $z^{\alpha}$ are holomorphic on $\Omega$, and they satisfy the differential relations 
$\frac{d}{dz} (\Lz) = \frac{1}{z}$ and $\frac{d}{dz}z^{\alpha} =\alpha z^{\alpha -1} $, respectively.

For a finite formal expression 
$$f(x) =\sum_{n \in \C, k \in \N} c_{n, k} \, x^{n} (\log x)^k \in \C_f\{x\}[\text{log}x]$$ and $z \in \Omega$, we define its evaluation at $z$ by
\[
f(z) = f(x)_{|x \leftarrow z} = \sum_{n, k} c_{n, k} \, z^{n} (\Lz)^k,
\]
obtained by replacing each occurrence of $x^n$ with $z^n$ and $(\log x)^k$ with $(\Lz)^k$.

This notion extends naturally to matrices. 
Let $A(x) = \big(a_{i,j}(x)\big)$ be an $m \times n$ matrix  whose entries lie in $\mathbb{C}_f\{x\}$.
For $\alpha \in \Omega$, we define its evaluation at $z$ as the matrix
\[
A(z)=A(x)_{| x \leftarrow z }=\big(a_{i,j}(x)_{| x \leftarrow z} \big)=(a_{i,j}(z)),
\]
obtained by evaluating each entry at $x=z$ as above.

The following lemma is a standard result for complex differential equations. For a proof, see \cite{W98}.

\begin{lemma}\label{diff-eq}
	Let $U$ be a simply connected domain in $\mathbb{C}$, and let $A(z)$ be an $n \times n$ matrix-valued holomorphic function on $U$. Then the space of all solutions to the differential equation $Y' = A(z)Y$ forms an $n$-dimensional vector space over $\mathbb{C}$. Here, $Y$ denotes a $\mathbb{C}^n$-valued holomorphic function on $U$.
\end{lemma}

The following property is referred to as the unique expansion set property in \cite{H25}. 
Here, we give an elementary proof for the finite sum case, which suffices for our purposes. 

\begin{lemma} \label{0=0}
	Let $\lambda_1, \lambda_2, \cdots, \lambda_n \in \mathbb{C}$ be distinct complex numbers, and $k_1, k_2, \cdots, k_m \in \mathbb{N}$ be distinct natural numbers.
	Let $U \subset \Omega$ be an open set. 
	Suppose the function
	\[
	F(z)=\sum_{i=1}^{n}\sum_{j=1}^{m}a_{i,j}z^{\lambda_i}(\Lz)^{k_j}, \ \ a_{i,j} \in \C, 
	\]
	satisfies $F(z)=0$ for all $z \in U$. Then  all coefficients $a_{i,j}=0$.
\end{lemma}

\begin{proof}
	
	Since the exponential function $e^z: \mathbb{C} \to \mathbb{C} \setminus \{0\}$ is continuous and surjective, there exists an open set $O \subset \mathbb{C}$ such that $\{e^z : z \in O\} \subset U$.  
	Therefore, the function  
	\[
	G(z) := F(e^z) = \sum_{i=1}^{n}\sum_{j=1}^{m} a_{i,j} z^{k_j} e^{\lambda_i z} = 0
	\]  
	for all $z \in O$.  
	
	Define the polynomials $f_i(z) = \sum_{j=1}^{m} a_{i,j} z^{k_j}$. Then  
	\[
	G(z) = \sum_{i=1}^{n} f_i(z) e^{\lambda_i z} = 0
	\]  
	for all $z \in O$.  
	To prove that $a_{i,j} = 0$ for all $i$ and $j$, it suffices to show that each polynomial $f_i(z)$ is identically zero on $O$.  
	
	Suppose, for contradiction, that there exists a function of the form
	\begin{align}\label{eqG=123}
    G(z) = \sum_{i=1}^{n} f_i(z) e^{\lambda_i z}
    \end{align}
    satisfying $G(z) = 0$ for all $z \in O$,
	where each $f_i(z)$ is a nonzero polynomial and $\lambda_1, \lambda_2, \cdots, \lambda_n$ 	
	are distinct complex numbers.
    By the method of minimal counterexample, we may assume that $G$ minimizes the value
	$$N_{G} := \big( \sum_{i=1}^n \deg f_i \big)+ n$$ 
	among all such functions.
	
	If necessary, multiply both sides of equation (\ref{eqG=123}) by $e^{-\lambda_1z}$. 
	Thus, without loss of generality, we may assume $\lambda_1=0$ and $\lambda_2, \cdots, \lambda_n$ 	
	are distinct, nonzero complex numbers. Differentiating $G(z)$ yields  
	$$0=G^{\prime}(z)=f_1^{\prime}(z)+ (\lambda_2 f_2(z) + f^{\prime}_2(z))e^{\lambda_2z} +\cdots + (\lambda_n f_n(z) + f^{\prime}_n(z))e^{\lambda_nz}. $$
	for all $z \in O$. Clearly, we have $N_{G^{\prime}(x)}<N_{G(x)}$. 
	By the minimality of $G$, we must have:
	$$\lambda_2 f_2(z) + f^{\prime}_2(z) = \cdots =\lambda_n f_n(z) + f^{\prime}_n(z)=0.$$
	However, this contradicts our assumption that each $f_i(z)$ is a nonzero polynomial and that $\lambda_2, \ldots, \lambda_n$ are nonzero complex numbers. Therefore, our initial assumption must be false, and consequently, if
	$$G(z) = \sum_{i=0}^{n} f_i(z) e^{\lambda_i z}=0$$
	for all $z \in O$, where each $f_i(z)$ is polynomial, then $f_i(z)=0$ for all $i$.
	
	This completes the proof.	
\end{proof}

\section{Preliminaries}\label{sec3}

Throughout this paper, we adopt the foundational framework of vertex operator algebras from \cite{FHL93, LL04}, 
with standard notation $(V, Y( -, x),  \1,  \w)$ for a vertex operator algebra. 

\subsection{Twisted Modules}

Let $(V, Y, \1, \omega)$ be a vertex operator algebra, and let $g$ be an automorphism of $V$ of finite order $T$. 
The $g$-eigenspace decomposition of $V$ is:
\[
V = \bigoplus_{r=0}^{T-1} V^{(g,r)}, \quad \text{where} \quad 
V^{(g,r)} = \left\{ v \in V \mid g v = e^{-2\pi\mathrm{i}r/T} v \right\}.
\]

\begin{definition}
A \textbf{ weak $g$-twisted $V$-module} $W$ is a vector space equipped with a linear map
	\[	Y_W(- , x) : V \to \operatorname{End}_{\C} (W)[[x^{\frac{1}{T}}, x^{-\frac{1}{T}}]], \quad \ \
v \mapsto Y_W(v, x)= \sum_{n \in \frac{1}{T}\Z} v_n x^{-n-1} \]
satisfying:	
\begin{enumerate}[{(1)}]
	\item $Y_W(\mathbf{1}, x) = \operatorname{id}_W$;
	
	\item For $u \in V^{(g,r)}$ and $w \in W$,  $Y_W(u, x)w \in x^{\frac{r}{T}} W((x))$;
	
	\item Jacobi identity: For $v \in V^r$, $u \in V$,
	\begin{align*}
		& x_0^{-1} \delta \left( \frac{x_1 - x_2}{x_0} \right) Y_W(v, x_1) Y_W(u, x_2) 
		- x_0^{-1} \delta \left( \frac{x_2 - x_1}{-x_0} \right) Y_W(u, x_2) Y_W(v, x_1) \\
		& \qquad = x_1^{-1}  \delta \left( \frac{x_2 + x_0}{x_1} \right) \left( \frac{x_2 + x_0}{x_1} \right)^{-\frac{r}{T}} Y_W(Y(v, x_0)u, x_2).
	\end{align*}
\end{enumerate}
\end{definition}

\begin{remark}
The definition of a weak $g$-twisted module given here corresponds to that of a $g^{-1}$-twisted module in \cite{DLM98}. 	
This convention, while not affecting the theory in essence, is purely for convenience.
\end{remark}

\begin{definition} \label{weight}
For a weak $g$-twisted $V$-module $W$ and $\lambda \in \mathbb{C}$, the $\lambda$-generalized eigenspace of $L(0)$ is defined as
\[
W_{(\lambda)} := \left\{ w \in W \ \middle| \ (L(0)-\lambda)^k w = 0 \ \text{ for some } k \in \mathbb{N} \right\}.
\]
A nonzero subspace $W_{(\lambda)}$ is called a {\bf weight space} of weight $\lambda$.  
If $w \in W_{(\lambda)} $ is nonzero, we write $\mathrm{wt}(w) = \lambda$.
For a homogeneous subspace $E \subseteq W$, we define its weight set as
\[
\mathrm{wt}(E) = \left\{ \lambda \in \mathrm{wt}(W) \mid E \cap W_{(\lambda)} \neq 0 \right\}.
\]

\end{definition}

\begin{definition}
	
	A {\bf generalized $g$-twisted $V$-module} is a weak $g$-twisted $V$-module $W$ admitting 
	a weight space decomposition $W = \bigoplus_{\lambda \in \mathbb{C}} W_{(\lambda)}$ 
	with the property that for every $\lambda \in \mathbb{C}$, the weight space $W_{(\lambda + \frac{1}{T} n)}$ is zero for all sufficiently negative integers $n$.
	
	If, in addition, each weight space $W_{(\lambda)}$ is finite-dimensional, 
	then $W$ is called a {\bf grading restricted generalized $g$-twisted $V$-module}.

    For a generalized $g$-twisted $V$-module $W$, let $W^{\prime}$ denote its restricted dual, regarded only as a vector space.
    
\end{definition}

Let $W$ be a generalized $g$-twisted $V$-module. For $\lambda \in \mathbb{C}$, define
$$
W[\lambda] = \bigoplus_{h \in \lambda + \frac{1}{T}\mathbb{Z}} W_{(h)}.
$$
Then $W[\lambda]$ is a generalized $g$-twisted $V$-submodule of $W$, and we have the direct sum decomposition
$$
W = \bigoplus_{\lambda \in \mathbb{C} / \frac{1}{T}\mathbb{Z}} W[\lambda].
$$
Consequently, if $W$ is indecomposable, there exists some $\lambda \in \mathbb{C}$ such that
$$
W = \bigoplus_{n \in \N} W_{(\lambda + \frac{n}{T})}.
$$


Let $W$ be a generalized $g$-twisted $V$-module. 
The semisimple part of the operator $L(0)$ is the linear map $L(0)_s : W \to W$ defined by
\[
 L(0)_s (w) =\lambda w  \ \ \text{for each } \ \  w \in W_{(\lambda)}.
\]
The nilpotent part of $L(0)$ is then defined as $L(0)_n = L(0) - L(0)_s$.

\vspace{0.5cm} 

The following lemma is well-known; for a proof, see \cite{M21}.

\begin{lemma}\label{lem:s-n}

Let $W$ be a generalized $g$-twisted $V$-module. The following statements hold:

\begin{enumerate}[{(1)}]
    \item Let $v \in V$ and $w \in W$ be homogeneous elements, and let $n \in \frac{1}{T}\mathbb{Z}$. Then $v_n w$ is homogeneous of weight
    \[
    \mathrm{wt}(v) - n - 1 + \mathrm{wt}(w).
    \]

    \item The linear map $L(0)_n \colon W \to W$ is an endomorphism of generalized $g$-twisted $V$-modules. In particular, if $W$ is finitely generated, then there exists $K \in \mathbb{N}$ such that
    \[
    (L(0)_n)^K W = 0.
    \]
\end{enumerate}
 
\end{lemma}
 
\begin{definition}

For a $g$-twisted $V$-module $W$, define the subspace $C_1(W)$ by
\begin{align*}
	C_1(W): = \ & \operatorname{span}_{\mathbb{C}} \left\{ v_{-1} w \mid v \in V^0,\ \operatorname{wt} v > 0,\ w \in W \right\} \\
	& + \operatorname{span}_{\mathbb{C}} \left\{ v_{-1 - r/T} w \mid v \in V^{(g, r)},\ r = 1, \cdots, T-1,\ \operatorname{wt} v \ge 0,\ w \in W \right\}.
\end{align*}

A $g$-twisted $V$-module $W$ is said to be {\bf $C_{1}$-cofinite} if the quotient space $W/C_{1}(W)$ is finite-dimensional.
\end{definition}

\begin{remark}\label{rmk:cofiniteness}
\begin{enumerate}[{(1)}]
    \item When $g = \id_V$, the definition of $C_1$-cofiniteness given above recovers the usual one (see, for example, \cite{M21, H05a}).

    \item Let $W$ be a nonzero generalized $g$-twisted $V$-module. Then $W$ admits a weight $\lambda_0$ for which none of the elements $\lambda_0 - n/T$ (with $n \in \mathbb{N}$) are weights of $W$. This implies, by the definition of $C_1(W)$, that $W_{(\lambda_0)} \cap C_1(W) = { 0 }$. Consequently, the quotient module $W / C_1(W)$ is nonzero.
    
    \item Let $W$ and $U$ be $g$-twisted $V$-modules,  Then there is a vector space isomorphism  
          \[
          (W \oplus U)/C_1(W \oplus U) \cong W/C_1(W) \oplus U/C_1(U).
          \]
          Consequently, if $W$ is a $C_1$-cofinite generalized $g$-twisted $V$-module, it must decompose into a direct sum of finitely many indecomposable submodules.

\end{enumerate}
\end{remark}

The proof of the following property parallels that for the untwisted case in \cite{M21}.

\begin{proposition}\label{c1-fg}

The following statement holds:
	\begin{enumerate}[{(1)}]

		\item Let $\lambda_1, \dots, \lambda_s \in \mathbb{C}$, $m \in \mathbb{N}$, and let $W$ be a $C_1$-cofinite generalized $g$-twisted $V$-module. Suppose $E$ is a homogeneous complementary subspace of $C_1(W)$ in $W$ such that
		\[
		\dim E \leq m \quad \text{and} \quad E \subset W_{(\lambda_1)} \oplus \cdots \oplus W_{(\lambda_s)}.
		\]
		Then for any $\lambda \in \mathbb{C}$, there exists a natural number $d_{\lambda}$---depending only on $\lambda_1, \dots, \lambda_s$, $m$, and the vertex operator algebra $V$, but independent of $W$---such that
		\[
		\dim W_{(\lambda)} \leq d_{\lambda}.
		\]

        \item Let $W$ be a $C_1$-cofinite generalized $g$-twisted $V$-module. Then $W$ is a grading restricted generalized $g$-twisted $V$-module and is finitely generated with at most $\dim (W/C_1(W))$ generators.
	\end{enumerate}
\end{proposition}

\begin{proof}

 (1) For any $\lambda \in \text{wt}  (W)$, let $\lambda_0$ be the weight of $W$ such that $\lambda-\lambda_0\in \frac{1}{T}\N$ is maximal. Denote this maximal value by $p$, 
and set
\[
q_{\lambda}=  \dim (\bigoplus_{0 \leq i \leq p } V_{(i)}).
\]
Since $W$ is $C_1$-cofinite, an induction argument then shows that
\[
W_{(\lambda)} \subset \text{Span}_{\C} \big\{v^1_{-1-r_1/T} \cdots v^n_{-1-r_n/T}w \ \ |
\ \  0 \leq \wt v^i \leq p, \ n \leq pT, \ w\in E   \big\}.
\]
It follows that $\dim (W_{(\lambda)} ) \leq (q_{\lambda})^{pT}m$,  which is independent of the twisted module $W$.
Let $d_{\lambda} = (q_{\lambda})^{pT}m$, as required. This completes the proof.

(2) This is a direct consequence of the proof of (1).
\end{proof}

\subsection{Twisted logarithmic Intertwining Operator and Fusion products}\label{sec3.2}

Let $g_1, g_2, g_3$  be three  mutually  commuting automorphisms of a vertex operator algebra $V$, 
of finite orders $T_1, T_2,  T_3$, respectively.
Since $g_1$ and $g_2$ commute, we have the following common eigenspace decomposition: 
\begin{equation*}
	V=\bigoplus_{0\leq r_1 < T_1, \  0 \leq r_2 < T_2}V^{(j_1, j_2)},
\end{equation*}
where 
\begin{equation*}\label{eq:common-eigenvector-space}
	V^{(r_1, r_2)}=\big\{ v \in V \mid  g_kv=e^{-2\pi \i r_k/T_k}v, \  k=1, 2 \big \}.
\end{equation*}

\begin{definition}\label{def:Intertwining}
Let $W^{k}$ be a $g_{k}$-twisted $V$-module for $k=1, 2, 3$.
A twisted logarithmic intertwining operator of type $\binom{W^3}{W^1 \ W^2}$ is a linear map
\begin{align*}
    \mathcal{Y}(- , x): & \ W^1 \otimes W^2  \  \to  \ W^3[\text{log}x]\{x\} \\
	&w^1 \otimes w^2 \ \to \  \Y(w^1, x)w^2= \sum_{n \in \C}\sum_{ k \in \N}(w^1)^\Y_{n;k}w^2x^{-n-1}(\text{log}x)^k \in  W^3[\text{log}x]\{x\}
\end{align*}
 satisfying the following conditions:

\begin{itemize}
\item \textbf{Lower truncation}: For any $w^1 \in W^1, w^2 \in W^2$, and $h \in \C,$ 
$$w^1_{h+n;k}w^2=0  \ \ \text{for} \ n \in \N \ \ \text{sufficiently large, independently of} \ k; $$
		
\item \textbf{Twisted Jacobi identity}: For $v \in V^{(r_1,r_2)}$, $w^1 \in W^1$, $w_2\in W^2$, and $0\leq r_1, r_2\leq T-1$,
    \begin{align} \label{Twisted Jacobi}
		&x_0^{-1} \delta(\frac{x_1-x_2}{x_0})(\frac{x_1-x_2}{x_0})^{-\frac{r_1}{T_1}}Y_{W^3}(v,x_1)\Y(w^1,x_2)w^2 \notag  \\
		-&x_0^{-1}\delta(\frac{-x_2+x_1}{x_0})(\frac{-x_2+x_1}{x_0})^{-\frac{r_1}{T_1}}\Y(w^1,x_2)Y_{W^2}(v,x_1)w^2 \notag  \\
		=&x_1^{-1}\delta(\frac{x_2+x_0}{x_1})(\frac{x_2+x_0}{x_1})^{-\frac{r_2}{T_2}}\Y(Y_{W^1}(v,x_0)w^1,x_2)w^2.  
	\end{align} 	  
	Here and below, expressions of the form $(-x)^\alpha$ are interpreted as $e^{\pi \i \alpha } x^\alpha$;
\item \textbf{$L(-1)$-derivative property}: For $w^1 \in W^1$, 
\begin{equation} \label{L(-1)-der}
			\Y(L(-1)w^1, x)=\frac{d}{d x}\Y(w^1, x).    
		\end{equation}
\end{itemize}

We simply denote $(w^1)^\Y_{n;k}w^2$ by $w^1_{n;k}w^2$  when there is no ambiguity.

The twisted logarithmic intertwining operators of a fixed type $\binom{W^3}{W^1 \ W^2}$ form a vector space, denoted by $\mathcal{I} \binom{W^3}{W^1 \ W^2}$. The dimension of this space is called the {\bf fusion rule} for $W^1$, $W^2$ and $W^3$, and is denoted by $N(W^1, W^2; W^3)$.

\end{definition}

\begin{remark}
Our definition of twisted logarithmic intertwining operators generalizes that of twisted intertwining operators in \cite{X95, DLM96, LS23, DLXY24, DXY25} and that of logarithmic intertwining operators in \cite{HLZ10}. 
We also note that a different approach, based on the language of convergence, 
has been used in \cite{DH25} to define such operators for non-commuting automorphisms.

\end{remark}

The following lemma, which follows from a direct computation, will be used in subsequent arguments.

\begin{lemma}
Following the notation of Definition \ref{def:Intertwining}, let $v \in V^{(r_1,r_2)}$. Then:

	\begin{enumerate}[{(1)}]
    \item  For $n, m \in \Z$, applying the residue operator 
$\mathrm{Res}_{x_0, x_1} x_0^{n-r_1/T_1} x_1^{m-r_2/T_2}$ to both sides of the twisted Jacobi identity yields:
	\begin{equation}\label{Operator-eq}
		\begin{aligned}
			& \sum_{i \geq 0 } (-1)^i \binom{n-r_1/T_1}{i} x_2^i v_{n-r_1/T_1+m-r_2/T_2-i} \Y(w^1, x_2)w^2 \\
			& - \sum_{i \geq 0 } \binom{n-r_1/T_1}{i} e^{\pi \i (n-r_1/T_1-i)}x_2^{n-r_1/T_1-i} \Y(w^1, x_2)v_{m-r_2/T_2+i}w^2\\
			&= \sum_{i  \geq 0 } \binom{m-r_2/T_2}{i} x_2^{m-r_2/T_2-i} \Y(v_{n-r_1/T_1+i}w^1, x_2)w^2. 
		\end{aligned}
	\end{equation}

\item 
Applying the residue operation $\mathrm{Res}_{x_2} x_2^h$ to both sides of \eqref{Operator-eq} 
and extracting the coefficient of $(\log x)^k$ yields the {\bf twisted Borcherds identity}:

\begin{equation}\label{twisted-B-id}
		\begin{aligned}
			& \sum_{i \geq 0 } (-1)^i \binom{n-r_1/T_1}{i}  v_{n-r_1/T_1+m-r_2/T_2-i} w^1_{h+i;k}w^2 \\
			& - \sum_{i \geq 0 } \binom{n-r_1/T_1}{i} e^{\pi \i (n-r_1/T_1-i)} w^1_{h+ n-r_1/T_1-i;k} v_{m-r_2/T_2+i}w^2\\
			&= \sum_{i  \geq 0 } \binom{m-r_2/T_2}{i} \big( v_{n-r_1/T_1+i}w^1 \big)_{h+m-r_2/T_2-i}w^2. 
		\end{aligned}
	\end{equation}

    \end{enumerate}

\end{lemma}

When the operator $\Y$ involves no $\text{log}x$ terms, the following property was established in \cite{X95}.

\begin{proposition} \label{g192=93}
We work in the setting of Definition \ref{def:Intertwining}. Assume that $Y_{W^3}(v, x)$ is nonzero for all $v \in V$, and that the logarithmic intertwining operator $\Y$ is nonzero. (Note that this first assumption holds, in particular, when $V$ is a simple vertex operator algebra.) Under these hypotheses, we have $g_3 = g_1 g_2$.
\end{proposition}

\begin{proof}
Take an arbitrary $v \in V^{(r_1, r_2)}$. From comparing the exponents of $x_1$ in the Jacobi identity, we deduce
\[
Y_{W^3}(v, x_1) \in x^{\frac{r_1}{T_1} + \frac{r_2}{T_2}} \operatorname{End}(W^3)((x_1)).
\]
Moreover, since $Y_{W^3}(u, x_3)$ is nonzero for every $u \in V$, it must be that $v \in V^{(g_3, r_3)}$, where the $r_3$ satisfies the condition
\[
\frac{r_3}{T_3} - \left( \frac{r_1}{T_1} + \frac{r_2}{T_2} \right) \in \mathbb{Z}.
\]
We therefore compute:
\[
g_3 v = e^{-2\pi i r_3 / T_3} v = e^{-2\pi i (r_1 / T_1 + r_2 / T_2)} v = g_1 g_2 v.
\]
This establishes the identity $g_3 = g_1 g_2$, and the proof is complete.
\end{proof}

\begin{definition}
Let $W^{k}$ be a $g_{k}$-twisted $V$-module for $k=1, 2, 3$.
Let $\Y$ be a twisted logarithmic intertwining operator of type $\binom{W^3}{W^1, W^2}$. Its image, $ \text{Im} \Y \subset W^3$, is defined as the subspace spanned by all coefficients in the series expansion of $\Y(w^1, x)w^2$ for all $w^1 \in W^1$ and $w^2 \in W^2$.
If $\text{Im}\Y =W^3$, then $\Y$ is called a surjective twisted logarithmic intertwining operator.
\end{definition}

\begin{proposition} \label{ImY-1}
The subspace $\text{Im} \Y$ is a $g_3$-twisted $V$-submodule of $W^3$.
\end{proposition}

\begin{proof}

 Let $w^1 \in W^1$, $w^2 \in W^2$, $l \in \mathbb{C}$, and $k \in \mathbb{N}$.
 By the lower truncation condition, there exists $n_0 \in \mathbb{N}$ such that
 \[
 w^1_{l+n_0+i;\,k}w^2 = 0 \quad \text{for all } i \in \mathbb{N}.
 \]
 
Let $v \in V$. We prove by induction on $d \in \N$ that for any $  q \in \frac{1}{T_3}\Z$, 
\begin{align*} \label{eq-123}
	v_{q}(w^1_{l+n_0-d;\,k}w^2) \in \text{Im} \Y.
\end{align*}

The base case $d = 0$ follows directly from the assumption that $w^1_{l+n_0;\,k}w^2 = 0$.

Now assume the statement holds for all $d < s$, where $s \geq 1$. 
That is, for any $q \in \frac{1}{T_3}\mathbb{Z}$ and $i \geq 1$:
\begin{align*} 
	v_{q}(w^1_{l+n_0-(s-i);\,k}w^2) \in \operatorname{Im} \mathcal{Y}
\end{align*}
 
Without loss of generality, assume $v \in V^{(r_1, r_2)}$ for some $r_1, r_2$, and that
\[
q = n - r_1 / T_1 + m - r_2 / T_2
\]
for some $m, n \in \mathbb{N}$. 
Setting $h=l+n_0-s$ in the twisted Borcherds identity \ref{twisted-B-id} and rearranging terms, we obtain

 \begin{align*}
 	&v_q (w^1_{l+n_0-s;\,k}w^2) =v_{n-r_1/T_1+m-r_2/T_2}(w^1_{l+n_0-s;\,k}w^2) \\
 	&\  = \sum_{i \geq 0} \binom{m-r_2/T_2}{i} (v_{n-r_1/T_1+i}w^1)_{l+n_0-s+m-r_2/T_2-i;\,k} w^2 \\
 	&\  + \sum_{i \geq 0} \binom{n-r_1/T_1}{i} e^{\pi \i (n-r_1/T_1-i)}  w^1_{n-r_1/T_1 + l+n_0-s-i} v_{m-r_2/T_2+i}w^2 \\
 	&\  - \sum_{i \geq 1} (-1)^i \binom{n-r_1/T_1}{i} v_{n-r_1/T_1+m-r_2/T_2-i} w^1_{l+n_0-(s-i)} w^2.
 \end{align*}
The first and second summations on the right-hand side clearly belong to $\text{Im} \Y$.
By the induction hypothesis, the third summation also lies in $\text{Im} \Y$.
Hence, the entire expression is contained in $\text{Im} \Y$.
This completes the induction step, and we conclude that
 \[
 v_{q}(w^1_{l+n_0+1-d;\,k}w^2) \in \text{Im} \Y \quad \text{for all } d \in \mathbb{N}.
 \]
 
Therefore $\text{Im} \Y$ is a weak $g_3$-twisted $V$-submodule of $W^3$.
The proof is complete.

\end{proof}

The proof of the following lemma is entirely analogous to that in the untwisted case given in \cite{HLZ10, M21}.
For completeness, we present it here.

\begin{lemma} \label{wt-id}
Following the notation of Definition~\ref{def:Intertwining}, suppose that $W^{k}$ is a generalized $g_{k}$-twisted $V$-module for $k=1, 2, 3$. We then state the following result.

\begin{enumerate}[{(1)}]

\item \label{wt-id:wt}
Let $w^1 \in W^1$ and $w^2 \in W^2$ be homogeneous elements, and let $h \in \mathbb{C}$. 
Then for every $k \in \mathbb{N}$, the element $w^1_{h;k} w^2$ is homogeneous of weight $\mathrm{wt}(w^1) + \mathrm{wt}(w^2) - h - 1$.

\item Let $d \in \N$ such that $\big(L(0)_n\big)^d W^i=0$ for $i=1,2,3$. Then there exists a $K \in \N$ such that 
\[
\Y(w^1, x)w^2 \in W^3\{x\}[\text{log}x]_{\leq K}
\]
for all $w^1 \in W^1$ and $w^2 \in W^2$.
\end{enumerate}

\end{lemma}

\begin{proof}

Setting $v = \omega \in V^{(0,0)}$ in the Jacobi identity and applying the residue operator $\mathrm{Res}_{x_0, x_1} x_1 $ to both sides of the twisted Jacobi identity, we obtain:
\begin{align*} 
&L(0)\Y(w^1, x_2)w^2 -\Y(w^1,x_2)(L(0)w^2) = x_2\Y(L(-1)w^1, x_2)w^2+ \Y(L(0)w^1, x_2)w^2 \\
& =x_2 \frac{d}{dx_2} \Y(w^1, x_2)w^2 + \Y(L(0)w^1, x_2)w^2
\end{align*}
Comparing the coefficients of  $x_2^{-h-1} (\text{log}x_2)^k$ yields
\begin{align*}
& L(0)(w^1_{h;k}w^2)-w^1_{h;k}(L(0)w^2) \\
& \ \  =(-h-1)w^1_{h;k}w^2+(k+1)w^1_{h;k+1}w^2+ (L(0)w^1)_{h;k}w^2.
\end{align*}
Consequently,
\begin{equation} \label{wt-eq}
\begin{aligned} 
&\big(L(0) - (\wt w^1 + \wt w^2 -h-1)\big) (w^1_{h;k}w^2) \\
\ & = \big((L(0) -\wt w^1 ) w^1 \big)_{h; k}w^2 + w^1_{h; k}\big((L(0)- \wt w^2)w^2\big) +(k+1)w^1_{h ;k+1}w^2.
\end{aligned}
\end{equation}

\noindent\textbf{(1)}
For $p \in  \N$,  applying  $\big(L(0) - (\wt w^1 + \wt w^2 -h-1) \big)^p$ to 
$w^1_{h;k}w^2$ and using mathematical induction, we can prove that
\begin{equation}
\begin{aligned} \label{wt-eq12}
&\big(L(0) - (\wt  w^1  + \wt w^2 -h-1) \big)^p (w^1_{h;k}w^2 ) \\
\ & =\sum_{i,j,s \in \N, i+j+s=p} \lambda_{i,j,s}\big((L(0) -\wt w^1 )^i w^1 \big)_{h; k+s} \big((L(0)- \wt w^2)^jw^2 \big),
\end{aligned}
\end{equation}
where $\lambda_{i,j,k} \in \C$.

By the definition of a generalized $g$-twisted module, together with the lower-truncation property of $Y(w^1, x)w^2$, there exists a $K \in \mathbb{N}$ such that
\[
(L(0)- \wt w^1)^K w^2=0, \ \ (L(0)- \wt w^2)^K w^2=0,
\]
and moreover, $w^1_{h;k}w^2=0 $  for any $k >K$. 
It follows from Equation \ref{wt-eq12} that 
\[
\left(L(0)- (\wt w^1 + \wt w^2 -h-1)\right)^{3K+1}(w^1_{h;k} w^2)=0.
\]
for any $k \in N$.
This completes the proof of (1).

\noindent\textbf{(2)}. Rearranging Equation \ref{wt-eq}, we obtain 
\[
(k+1)w^1_{h; k+1}w^2= (L(0)_nw^1)_{h;k}w^2 + w^1_{h;k} (L(0)_nw^2) -L(0)_n(w^1_{h; k}w^2) 
\]
It follows that for any $s= 0, 1, \cdots,   k+1$ we have  
\begin{equation} \label{wt-eq123}
\begin{aligned}
w^1_{h;k+1}w^2=\sum_{a,b,c,s \in \N, a+b+c=s } \lambda_{a,b,c}^s L(0)_n^a \left( \big( L(0)_n^bw^1 \big)_{h, k+1-s} \big( L(0)_n^c w^2 \big) \right),
\end{aligned}
\end{equation}
where $\lambda_{a,b,c}^s \in \C$.
Taking $k = 3d$ and $s=3d+1$ in the above identity, we obtain 
\[w^1_{h; 3d+1}w^2 = 0\] for all $w^1 \in W^1$ and $w^2 \in W^2$. 
Let $K=3d+1$. Then we have  
\[
\mathcal{Y}(w^1, x) w^2 \in W^3\{x\}[\log x]_{\leq K}
\]
for any $w^1 \in W^1$ and $w^2 \in W^2$. This completes the proof of (2).
\end{proof}

\begin{remark}
Let $W^{k}$ be a generalized $g_{k}$-twisted $V$-module for $k=1, 2, 3$.
Let $\Y$ be a twisted logarithmic intertwining operator of type $\binom{W^3}{W^1, W^2}$.
Let $\theta \in (W^3)^{\prime}, w^1 \in W^1$ and $w^2 \in W^2$. It follows from Lemma \ref{wt-id}(\ref{wt-id:wt})
that
\[
\big< \theta, \ \mathcal{Y}(w^1, x)w^2 \big> =\sum_{n \in \C, k \in \N} \big< \theta, \ w^1_{n;k}w^2  \big> x^{-n-1}(\text{log}x)^k
\]
is a finite sum. Hence
\[
\big< \theta, \ \mathcal{Y}(w^1, z)w^2 \big> = \big< \theta, \ \mathcal{Y}(w^1, x)w^2 \big>_{|x \leftarrow z} 
\]
is a holomorphic function on $\Omega$. This fact will be important in what follows.

\end{remark}

\begin{definition}
Let $W^1$ and $W^2$ be generalized $g_1$-twisted and $g_2$-twisted 
$V$-modules, respectively. A \textbf{fusion product} of $W^1$ and $W^2$ 
is a pair $(W^1\boxtimes W^2, \mathcal{Y}^{\boxtimes})$, where $W^1\boxtimes W^2$ is a generalized $g_3$-twisted 
$V$-module and $\mathcal{Y}^{\boxtimes}$ is a twisted logarithmic intertwining operator of type 
$\binom{W^1\boxtimes W^2}{W^1 \ W^2}$, satisfying the following universal property: 
For any generalized $g_3$-twisted $V$-module $W$ and any twisted logarithmic 
intertwining operator $\mathcal{Y}$ of type $\binom{W}{W^1 \ W^2}$, there 
exists a unique generalized $g_3$-twisted $V$-module homomorphism $\phi: W^1\boxtimes W^2 \to W$ such that
\[\mathcal{Y} = \phi \circ \mathcal{Y}^{\boxtimes}.
\]
\end{definition}

\begin{remark}
The treatment of the fusion product in this work is confined to the setting of generalized modules.
\end{remark}

\begin{proposition}
Assume that the fusion product $W^1 \boxtimes W^2$ exists. Then $\Y^{\boxtimes}$ is a surjective twisted logarithmic intertwining operator.
\end{proposition}
\begin{proof}
The proof is routine. By Proposition~\ref{ImY-1}, the image $\text{Im}\,\Y^{\boxtimes}$ is a generalized $g_3$-twisted $V$-submodule of the fusion product $W^1 \boxtimes W^2$. Consequently, the quotient $W := (W^1 \boxtimes W^2) / \text{Im}\,\Y^{\boxtimes}$ is a generalized $g_3$-twisted $V$-module. Let $\pi: W^1 \boxtimes W^2 \to W$ be the canonical quotient map.

We then have $\pi \circ \Y^{\boxtimes} \in \mathcal{I} \binom{W}{W^1 \, W^2}$. Noting that $\pi \circ \Y^{\boxtimes} = 0$, the universal property of the fusion product implies that $\pi = 0$. Therefore, $\text{Im}\,\Y^{\boxtimes} = W^1 \boxtimes W^2$, which proves that $\Y^{\boxtimes}$ is surjective. This completes the proof.
\end{proof}

\section{ \texorpdfstring{$C{_1}$}{C1}-cofiniteness of Fusion Product}\label{sec4}

Let $g_1, g_2$ be two commuting automorphisms of a vertex operator algebra $V$, of finite orders $T_1, T_2$, respectively. Set $g_3=g_1g_2$. Let $T_3$ be the order of $g_3$.
For $k = 1, 2, 3$, let $W^k$ be a generalized $g_k$-twisted $V$-module, and let $E_k$ be a homogeneous complementary subspace of $C_1(W^k)$ in $W^k$, so that
\[
W^k = E_k \oplus C_1(W^k), \quad \text{for } k = 1, 2, 3.
\]
For the rest of this section, we assume that both $W^1$ and $W^2$ are $C_1$-cofinite.

Define
\[
E_3^\circ = \left\{ f \in (W^3)' \ \middle| \ f(C_1(W^3)) = 0 \right\} \subset (W^3)'.
\]
where $(W^3)'$ is the restricted dual of $W^3$.

\begin{lemma} \label{man-lemma}

Let \(\mathcal{Y}\) be a twisted logarithmic intertwining operator of type \(\binom{W^3}{W^1\,W^2}\), 
\(\theta \in E_3^\circ\), and let \(w^1 \in W^1\), \(w^2 \in W^2\) be homogeneous elements. 
Then there exist finitely many elements \(p_i \in E_1\), \(q_j \in E_2\) and polynomials
\[
f_{i,j}(x) \in \mathbb{C}\big[x^{\pm 1/T_1}, x^{\pm 1/T_2}\big],
\] 
depending only on \(w^1\) and \(w^2\) (and independent of \(\mathcal{Y}\) and \(\theta\)), such that
\[
\big\langle \theta, \mathcal{Y}(w^1, x)w^2 \big\rangle = \sum_{i,j} \big\langle \theta, \mathcal{Y}(p_i, x)q_j \big\rangle f_{i,j}(x).
\]

\end{lemma}

\begin{proof}

Given the lower truncation of weights and the homogeneity of $w^1$ and $w^2$, 
there exist weights $\lambda_1 \in \wt(W^1)$ and $\lambda_2 \in \wt(W^2)$ such that $\wt(w^1) - \lambda_1$ and $\wt(w^2) - \lambda_2$ are maximal in $\frac{1}{T_1}\mathbb{N}$ and $\frac{1}{T_2}\mathbb{N}$, respectively.

We prove the lemma by induction on $\wt(w^1)-\lambda_1 + \wt(w^2)-\lambda_2$. The base case is immediate from the inclusions $W^1_{(\lambda_1)} \subset E_1$ and $W^2_{(\lambda_2)} \subset E_2$, which hold by the choice of $\lambda_1$ and $\lambda_2$. For the inductive step, we assume without loss of generality that $w^1 \in C_1(W^1)$ or $w^2 \in C_1(W^2)$.


\vspace{0.5\baselineskip}

{\bf Case 1: }	$w^1 \in C_1(W^1)$.
By the definition of $C_1(W^1)$ and the linearity of the form $\langle \theta, \mathcal{Y}(-, x)w^2 \rangle$ 
	in its first argument, we may assume without loss of generality that $w^1 = v_{-1 - r_1/T_1}p$ 
	for some homogeneous $v \in V^{(r_1, r_2)}$ and $p \in W^1$. 
	Then
	\[
	\operatorname{wt} w^1 + \operatorname{wt} w^2 = \operatorname{wt} v + \frac{r_1}{T_1} + \operatorname{wt} p + \operatorname{wt} w^2.
	\]
	Setting $m = -1, n = 0$, and replacing $x_2$  with $x$ in identity \ref{Operator-eq}, we obtain:
	\begin{equation}
		\begin{aligned}
			&\big< \theta, \mathcal{Y}(w^1, x)w^2 \big> = \big< \theta, \mathcal{Y}(v_{-1 - r_1/T_1}p, x)w^2 \big> \\
			& \ = \sum_{i \geq 0} (-1)^i \binom{-1 - r_1/T_1}{i} x^i 
			\big< \theta, v_{-1 - r_1/T_1 - r_2/T_2 - i} \mathcal{Y}(p, x)w^2 \big> \\
			&\ - \sum_{i \geq 0} \binom{-1 - r_1/T_1}{i} e^{\pi \i (-1 - r_1/T_1 - i)}x^{-1 - r_1/T_1 - i} 
			\big< \theta, \mathcal{Y}(p, x)v_{-r_2/T_2 + i}w^2 \big> \\
			&\ - \sum_{i \geq 1} \binom{-r_2/T_2}{i} x^{-r_2/T_2 - i} 
			\big< \theta, \mathcal{Y}(v_{-1 - r_1/T_1 + i}p, x)w^2 \big>.
		\end{aligned}
	\end{equation}
	
	Since $g_3=g_1g_2$, we have
	\[
	\frac{r_1}{T_1} + \frac{r_2}{T_2} = \frac{r_3}{T_3} \quad \text{or} \quad \frac{r_1}{T_1} + \frac{r_2}{T_2} = \frac{r_3}{T_3} + 1
	\] 
	for some $0 \leq r_3 < T_3$. Then for $i \geq 0$, we have
	\[
	v_{-1 - r_1/T_1 - r_2/T_2 - i} \mathcal{Y}(p, x) w^2 \in C_1(W^3)[\text{log}x]\{x\},
	\]
	and thus
	\[
	\big< \theta, v_{-1 - r_1/T_1 - r_2/T_2 - i} \mathcal{Y}(p, x) w^2 \big> = 0
	\]
	by the choice of $\theta$. Moreover, we observe that
	\[
	\operatorname{wt} p + \operatorname{wt}(v_{-r_2/T_2 + i} w^2) < \operatorname{wt} w^1 + \operatorname{wt} w^2 \quad \text{for any } i \geq 0,
	\]
	and
	\[
	\operatorname{wt}(v_{-1 - r_1/T_1 + i} p) + \operatorname{wt} w^2 < \operatorname{wt} w^1 + \operatorname{wt} w^2 \quad \text{for any } i \geq 1.
	\]
	
	By the induction hypothesis, there exist elements $p^1_{i,a}, p^2_{i,a} \in E_1$ and $q^1_{i,b}, q^2_{i,b} \in E_2$, 
	together with elements $f_{i;a,b}(x), g_{i;a,b}(x) \in \mathbb{C}[x^{\pm1/T_1}, x^{\pm 1/T_2}]$, 
	all independent of $\mathcal{Y}$ and $\theta$, such that
	\[
	\big< \theta, \mathcal{Y}(p, x)v_{-r_2/T_2 + i}w^2 \big> = \sum_{a, b} \big< \theta, \mathcal{Y}(p^1_{i,a}, x)q^1_{i,b} \big> f_{i;a,b}(x),
	\]
	and
	\[
	\big< \theta, \mathcal{Y}(v_{-1 - r_1/T_1 + i}p, x)w^2 \big> = \sum_{a,b} \big< \theta, \mathcal{Y}(p^2_{i,a}, x)q^2_{i,b} \big> g_{i;a,b}(x).
	\]
  Substituting these expressions back into the earlier formula, we see that $\langle \theta, \mathcal{Y}(w^1, x)w^2 \rangle$ can be written as a finite linear combination of terms of the form $\langle \theta, \mathcal{Y}(p_i, x)q_j \rangle f_{i,j}(x)$, where each $p_i \in E_1$, $q_j \in E_2$, and $f_{i,j}(x) \in \mathbb{C}[x^{\pm1/T_1}, x^{\pm 1/T_2}]$, and these elements depend only on $w^1$ and $w^2$, not on $\mathcal{Y}$ and $\theta$. This completes the inductive step in this case.

	

\vspace{0.5\baselineskip}

{\bf Case 2: }	$w^2 \in C_1(W^2)$.
	By symmetry, we may assume that $w^2 = v_{-1 - r_2/T_2}q$ for some homogeneous $v \in V^{(r_1, r_2)}$ and $q \in W^2$. 
	Setting $m = 0, n = -1$, and replacing $x_2$  with $x$ in identity \ref{Operator-eq}, we obtain:
	\begin{equation}
		\begin{aligned}
			&\big< \theta, \mathcal{Y}(w^1, x)w^2 \big> = \big< \theta, \mathcal{Y}(w^1, x)v_{-1 - r_2/T_2}q \big> \\
			& \ = \sum_{i \geq 0} (-1)^i \binom{-1 - r_2/T_2}{i} x^i 
			\big< \theta, v_{-1 - r_1/T_1 - r_2/T_2 - i} \mathcal{Y}(w^1, x)q \big> \\
			& \ - \sum_{i \geq 0} \binom{-1 - r_2/T_2}{i}e^{\pi \i (-1 - r_2/T_2 - i)} x^{-1 - r_2/T_2 - i} 
			\big< \theta, \mathcal{Y}(w^1, x)v_{-r_1/T_1 + i}q \big> \\
			&\ - \sum_{i \geq 1} \binom{-r_1/T_1}{i} x^{-r_1/T_1 - i} 
			\big< \theta, \mathcal{Y}(v_{-1 - r_2/T_2 + i}w^1, x)q \big>.
		\end{aligned}
	\end{equation}
	
	A similar argument shows that the lemma holds in this case as well.

This completes the proof.		
\end{proof}

\vspace{0.5\baselineskip}

Assume $\dim E_1 = s$ and $\dim E_2 = t$, and let $m = st$. Let $\{p^1, \dots, p^s\}$ be a basis for $E_1$ and $\{q^1, \dots, q^t\}$ a basis for $E_2$.
Let $\Y$ be a twisted logarithmic intertwining operator of type $\binom{W^3}{W^1 \ W^2}$, and $\theta \in E_3^\circ$.

For $1 \leq n \leq m$, write $n = (i-1)t + j$ where $1 \leq i \leq s$ and $1 \leq j \leq t$, and define
\[
Y_n(\Y, \theta; x) =  \big<  \theta, \Y(p^i, x)q^j \big> \in \mathbb{C}_f\{x\}[\log x],
\]
and let
\[
Y(\Y, \theta; x) = \big( Y_1(\Y, \theta; x), \cdots, Y_m(\Y, \theta; x) \big)^T \in \big( \mathbb{C}_f\{x\}[\log x] \big)^{\oplus m}.
\]
By Lemma \ref{man-lemma}, for any $1 \leq i \leq s$ and $1 \leq j \leq t$, there exist elements 
\[
f_{a,b}^{i,j}(x) \in \mathbb{C}[x^{\pm 1/T_1}, x^{\pm 1/T_2}], \ \ 1 \leq a \leq s, 1 \leq b \leq t,
\]
independent of $\Y$ and $\theta$, such that
\[
\frac{d}{dx} \big< \theta, \Y(p^i, x)q^j \big> = \big< \theta, \Y(L(-1)p^i, x)q^j \big> 
= \sum_{a,b} f_{a,b}^{i,j}(x) \big< \theta, \Y(p^a, x)q^b \big>.
\]
Consequently, there exists an $m \times m$ matrix $A(x)$ with entries in $\mathbb{C}[x^{\pm 1/T_1}, x^{\pm 1/T_2}]$, also independent of $\Y$ and $\theta$, such that
\[
\frac{d}{dx} Y(\Y, \theta; x) = A(x) Y(\Y, \theta; x).
\]

For $z \in \Omega$, define
\[
Y(\Y, \theta; z) =Y\big(\Y, \theta; x\big)_{|x \leftarrow z}.
\]
Then
\[
\frac{d}{dz} Y(\Y, \theta; z) = \big(\frac{d}{dx} Y\big(\Y, \theta; x\big)\big)_{|x \leftarrow z} = A(z) Y(\Y, \theta; z).
\]
Therefore, for any twisted logarithmic intertwining operator $\Y$ of type $\binom{W^3}{W^1 \ W^2}$ 
and any $\theta \in E_3^\circ$, the corresponding vector-valued function $Y(\Y, \theta; z)$ satisfies the linear differential equation
\[
\frac{d}{dz} Y( z) = A(z) Y( z)
\]
on the simply connected domain $\Omega$.

Let $H$ be the solution space of the linear differential equation
\[
\frac{d}{dz} Y(z) = A(z) Y(z)
\]
on the simply connected domain $\Omega$. 
By Lemma \ref{diff-eq}, we have $\dim H = m$. 
Moreover, $Y(\mathcal{Y}, \theta; z) \in H$ 
for every such $\mathcal{Y}$ and  such $\theta \in E_3^\circ$.

\begin{theorem}\label{cofin-th}
Suppose there exists a twisted logarithmic intertwining operator $\mathcal{Y}$ of type $\binom{W^3}{W^1 \ W^2}$ satisfying $\operatorname{Im} \mathcal{Y} = W^3$.
Then $\dim E_3^\circ \leq m$, and hence $W^3$ is $C_1$-cofinite. In particular, $W^1 \boxtimes W^2$ is a $C_1$-cofinite $g_3$-twisted generalized $V$-module.
\end{theorem}

\begin{proof}

	Define a linear map $f \colon E_3^\circ \to H$ by $f(\theta) = Y(\mathcal{Y}, \theta; z)$. 
	Since $\text{dim} H =m$, to prove the theorem, it suffices to show that $f$ is injective.
	
	Suppose $\theta \in \ker f$. Then 
	\[
	\langle \theta, \mathcal{Y}(p^i, z)q^j \rangle = 0 \  \ ( 1 \leq i \leq s, \  1 \leq j \leq t) 
	\]
    for all $ z \in \Omega$.
    By Lemma \ref{man-lemma}, for any $w^1 \in W^1$ and $w^2 \in W^2$, 
    there exist 
    \[
    f_{i,j}(x) \in \mathbb{C}[x^{\pm 1/T_1}, x^{\pm 1/T_2}] \ \ \text{for} \ \ 1 \leq i \leq s, 1 \leq j \leq t 
    \]
    such that
	\[
	\langle \theta, \mathcal{Y}(w^1, x)w^2 \rangle = \sum_{i,j} \langle \theta, \mathcal{Y}(p^i, x)q^j \rangle f_{i,j}(x).
	\]
	Consequently,
	\[
	\langle \theta, \mathcal{Y}(w^1, z)w^2 \rangle = \sum_{i,j} \langle \theta, \mathcal{Y}(p_i, z)q_j \rangle f_{i,j}(z) = 0
	\]
	for all $z \in \Omega$. This implies that the finite sum
	\[
	\sum_{n \in \mathbb{C},\, k \in \mathbb{N}} \langle \theta, w^1_{n;k} w^2 \rangle z^{-n-1} (\Lz)^k = 0
	\]
	for all $ z \in \Omega$. Applying Lemma \ref{0=0}, we conclude that
	\[
	\langle \theta, w^1_{n;k} w^2 \rangle = 0
	\]
	for all $n \in \mathbb{C}$ and $k \in \mathbb{N}$. Hence,
	\[
	\langle \theta, W^3 \rangle = \langle \theta, \operatorname{Im} \mathcal{Y} \rangle = 0,
	\]
	which implies $\theta = 0$. Therefore, $f$ is injective, completing the proof.
\end{proof}

\begin{theorem}\label{finite-dim}
Let $W^k$ (for $k=1,2$) be a $C_1$-cofinite generalized $g_k$-twisted $V$-module, and let $W^3$ be a generalized $g_3$-twisted $V$-module of finite composition length. Then the fusion rule $N(W^1, W^2; W^3)$ is finite-dimensional.
\end{theorem}

\begin{proof}

We first prove the theorem in the case where $W^3$ is a simple module.

Fix a nonzero element $\theta \in E_3^\circ$, whose existence is guaranteed by Remark~\ref{rmk:cofiniteness}(2), and define a linear map
\[
f \colon \mathcal{I} \binom{W^3}{W^1 \; W^2} \to H, \quad \text{by} \quad f(\mathcal{Y}) = Y(\mathcal{Y}, \theta; z).
\]
Suppose $\mathcal{Y} \in \ker f$. Then
\[
\langle \theta, \mathcal{Y}(p^i, z) q^j \rangle = 0 \quad \text{for all }  z \in \Omega,  \  1 \leq i \leq s, \  1\leq j \leq t.
\]
By the same argument as in the proof of Theorem~\ref{cofin-th}, we conclude that  
\[
\langle \theta, \operatorname{Im} \mathcal{Y} \rangle = 0.
\]
Since $W^3$ is simple and $\operatorname{Im} \mathcal{Y} \leq W^3$, it follows that $\operatorname{Im} \mathcal{Y} = 0$, and hence $\mathcal{Y} = 0$.  
Thus, $f$ is injective, and we obtain
\[
N(W^1, W^2; W^3) \leq \dim H = m,
\]
which shows that the fusion rule $N(W^1, W^2; W^3)$ is finite-dimensional when $W^3$ is simple.

Now suppose $W^3$ has composition length at least $2$. Consider the short exact sequence
\[
0 \to S \xrightarrow{j} W^3 \xrightarrow{\pi} M \to 0.
\]
Since the composition lengths of $S$ and $M$ are strictly less than that of $W^3$, the induction hypothesis implies that
\[
N(W^1, W^2; S) < \infty \quad \text{and} \quad N(W^1, W^2; M) < \infty.
\]
Define a linear map
\[
f \colon \mathcal{I} \binom{W^3}{W^1 \; W^2} \to \mathcal{I} \binom{M}{W^1 \; W^2} \quad \text{by} \quad f(\mathcal{Y}) = \pi \mathcal{Y}.
\]
If $\mathcal{Y} \in \ker f$, then $\pi(w^1_{n;k} w^2) = 0$ for all $w^1 \in W^1$, $w^2 \in W^2$, $n \in \mathbb{C}$, and $k \in \mathbb{N}$. So $\operatorname{Im} \mathcal{Y} \subset S$. Therefore, $\ker f \subset \mathcal{I} \binom{S}{W^1 \; W^2}$, and we obtain
\[
N(W^1, W^2; W^3) \leq N(W^1, W^2; S) + N(W^1, W^2; M) < \infty.
\]

This completes the proof.

\end{proof}

\begin{remark}

The finiteness of the fusion rule $N(W^1, W^2; W^3)$ was established in \cite{L99} for irreducible $g$-twisted $V$-modules, and later proven in \cite{H05a} under the broader assumption that $W^1$, $W^2$, and $(W^3)'$ are $C_1$-cofinite.
\end{remark}

\section{Constructions of the fusion product}\label{sec5}


Let $g_1, g_2$ be two commuting automorphisms of a vertex operator algebra $V$, of finite orders $T_1, T_2$, respectively. 
Define $g_3 = g_1 g_2$, and let $T_3$ denote the order of $g_3$. 
For $k=1,2$, let $W^k$ be a $C_1$-cofinite generalized $g_k$-twisted $V$-module 
with codimension 
\[
m_k = \dim\left(W^k/C_1(W^k)\right).
\]
Define $m = m_1 m_2$.
From now on, we fix $W^1, W^2$ and $m$.

Let 
\[
\Y_1 \in \I \binom{U}{W^1 \ W^2} \  \  \text{and}  \  \  \Y_2 \in \I \binom{W}{W^1 \ W^2} 
\]
be two surjective twisted logarithmic intertwining operators.
We say the pair $(U, \Y_1)$ and $(W, \Y_2)$ are equivalent, denoted by $(U, \Y_1) \cong (W, \Y_2)$, 
if there exists an isomorphism $f: U \to W$ of $g_3$-twisted $V$-modules such that $f\Y_1=\Y_2$.

Now, define the set
\[
\mathcal{F}_{(W^1,W^2)}= \left\{ (U, \mathcal{Y}) \ \middle| \  U \in  V\text{-}\mathbf{Mod}^{g_3} \ 
\text{and} \ \mathcal{Y}  \in \mathcal{I} \binom{U}{W^1 \ W^2} \ \text{is surjective} \right\} \Big/ \cong,
\]
where $V\text{-}\mathbf{Mod}^{g_3}$ denotes the category of generalized $g_3$-twisted $V$-modules.
This indeed constitutes a set, since all modules $U$ involved are finitely generated with at most $m$ generators (Theorem \ref{cofin-th} and Proposition \ref{c1-fg}).
Therefore, their isomorphism classes form a set.
Moreover, for any such module $U$, the collection of surjective intertwining operators $\mathcal{Y}$ also forms a set.
Note that if $(W, \Y) \in \mathcal{F}_{(W^1,W^2)}$, then $W$ is automatically a grading restricted generalized $g_3$-twisted $V$-module.

\subsection{General case} \label{sec6-1}

Let \( W \) be a generalized \( g \)-twisted \( V \)-module. The subspace \( C_1(W) \) is invariant under \( L(0)_s \), which therefore induces a natural action on the quotient \( W/C_1(W) \). Denote by \( \operatorname{wt}(W/C_1(W)) \) the eigenvalue set of \( L(0)_s \) on this quotient.
We note two facts:
\begin{enumerate}[{(1)}]
    \item If \( E \) is any homogeneous complement of \( C_1(W) \) in \( W \), then \( \operatorname{wt}(E) = \operatorname{wt}(W/C_1(W)) \);
    \item Let \( f: W \to U \) be a homomorphism of generalized \( g \)-twisted \( V \)-modules. The condition \( f(W_{(\lambda)}) \subset U_{(\lambda)} \) implies that \( f \) commutes with \( L(0)_s \). Consequently, if \( f \) is surjective, then we have
    $$\wt(U/C_1(U)) \subset \wt(W/C_1(W)).$$
\end{enumerate}





\begin{lemma} \label{leq-d}

There exist finitely many complex numbers $\lambda_1, \dots, \lambda_n $ and a positive integer $d \in \mathbb{N}$ such that for every $(W, \Y) \in \mathcal{F}_{(W^1,W^2)}$, the following hold:
\begin{enumerate}[{(1)}]
\item Every homogeneous complementary subspace $E_W$ of $C_1(W)$ in $W$ is contained in $\bigoplus_{i=1}^n W_{(\lambda_i)}$;
\item The direct sum $\bigoplus_{i=1}^n W_{(\lambda_i)}$ has dimension at most $d$.
\end{enumerate}

\end{lemma}

\begin{proof}

(1) By Theorem~\ref{cofin-th}, we have 
\[
\dim(W/C_1(W)) \leq m \ \ \ \text{for any} \ \ \ (W, \mathcal{Y}) \in \mathcal{F}_{(W^1,W^2)}.
\]
Hence, there exists $(U, \mathcal{Y}_U) \in \mathcal{F}_{(W^1,W^2)}$ such that $\dim(U/C_1(U))$ is maximal.  
Let
\[
\wt(U/C_1(U)) = \{\lambda_1, \dots, \lambda_n\}.
\]

Now, for any $(W, \mathcal{Y}) \in \mathcal{F}_{(W^1,W^2)}$, we relate $W$ to $U$ via the direct sum intertwining operator
\[
\mathcal{Y}_U \oplus \mathcal{Y} \in \mathcal{I} \binom{U \oplus W}{W^1\, W^2}.
\]
Let $p \colon U \oplus W \to U$ and $q \colon U \oplus W \to W$ be the canonical projections, and define
\[
G = \mathrm{Im}(\mathcal{Y}_U \oplus \mathcal{Y}), \quad \mathcal{Y}^\circ = \mathcal{Y}_U \oplus \mathcal{Y}.
\]
Then, $(G, \mathcal{Y}^\circ) \in \mathcal{F}_{(W^1,W^2)}$, and the restriction $p|_G \colon G \to U$ is an epimorphism.  
This induces a surjective linear map
\[
\widetilde{p} \colon G/C_1(G) \to U/C_1(U).
\]
By the maximality of $\dim U/C_1(U)$, the map $\widetilde{p}$ is an isomorphism.  
Consequently,
\[
\wt(G/C_1(G)) = \wt(U/C_1(U)) = \{\lambda_1, \dots, \lambda_n\}.
\]

On the other hand, $q|_G \colon G \to W$ is also an epimorphism. Therefore,
\[
\wt(W/C_1(W)) \subset \wt(G/C_1(G)) = \{\lambda_1, \dots, \lambda_n\}.
\]
It follows that for any homogeneous complementary subspace $E_W$ of $C_1(W)$ in $W$,
\[
\wt(E_W) \subset \{\lambda_1, \dots, \lambda_n\},
\]
and hence $E_W$ is contained in $\bigoplus_{i=1}^n W_{(\lambda_i)}$.
This completes the proof of (1).

(2)  By Lemma \ref{leq-d}(1) and Lemma~\ref{c1-fg}, there exist $d_{\lambda_1}, \dots, d_{\lambda_n} \in \mathbb{N}$ such that for any $(W, \mathcal{Y}) \in \mathcal{F}_{(W^1,W^2)}$ and each $i$,
\[
\dim W_{(\lambda_i)} \leq d_{\lambda_i}.
\]
Setting $d = d_{\lambda_1} + \cdots + d_{\lambda_n}$, we conclude that for all such $W$,
\[
\dim \bigoplus_{i=1}^n W_{(\lambda_i)} \leq d.
\]

This completes the proof of (2).

\end{proof}

From now on, we fix these $\lambda_1, \cdots, \lambda_n$ and $d$.

\begin{lemma} \label{U-Y-Lemma}

    For any $(W, \mathcal{Y}) \in \mathcal{F}_{(W^1,W^2)}$, the following statements hold:
    
    \begin{enumerate}[{(1)}]
        \item The set of weights satisfies
        \[
        \mathrm{wt}(W) \subset \bigcup_{i=1}^n \left( \lambda_i + \frac{1}{T_3} \mathbb{N} \right).
        \]
        
        \item There exists a $K \in \N$ such that for any $w^1 \in W^1$ and $w^2 \in W^2$,
        \[
        \mathcal{Y}(w^1, x) w^2 \in U\{x\}[\log x]_{\leq K}.
        \]
    \end{enumerate}
		
\end{lemma}

\begin{proof}
(1) Let $E_W$ be a homogeneous complement of $C_1(W)$ in $W$. 
Then $W$ is spanned by elements of the form
\[
v^1_{-1-r_1/T_3}v^2_{-1-r_2/T_3} \cdots v^s_{-1-r_s/T_3}e, \quad 
s \in \mathbb{N}, \ v^i \in V, \ \operatorname{wt}(v^i) \geq 0, \ e \in E_W.
\]
The result now follows immediately from Lemma~\ref{leq-d}(1).

(2) By Lemma \ref{leq-d}, we know that $W$ is generated by $\bigoplus_{i=1}^n W_{(\lambda_i)}$, and we have
\[
\left(L(0)_n\right)^d \big(\bigoplus_{i=1}^n W_{(\lambda_i)} \big) =0.
\]
Consequently, by Lemma \ref{lem:s-n}, we have
\[
\left(L(0)_n\right)^d W=0 \ \ \text{ for any} \ \ (W, \mathcal{Y}) \in \mathcal{F}_{(W^1,W^2)}.
\]
Since both $W^1$ and $W^2$ are $C_1$-cofinite, they are finitely generated. 
Therefore, there exists $p \in \mathbb{N}$ such that
\[
(L(0)_n)^p W^i=0 \ \ \text{for} \ \  i=1, 2.
\] 
The result now follows immediately from Lemma \ref{wt-id}(2).

\end{proof}

We now proceed to construct the fusion product of $W^1$ and $W^2$.
Consider the direct product space
\[
\mathcal{S} = \prod_{(U, \mathcal{Y}) \in \mathcal{F}(W^1,W^2)} U.
\]
In this space, we define the vertex operator $Y_{\mathcal{S}}(-, x)$ by taking the direct product of the vertex operators on each component:
\[
Y_{\mathcal{S}}(v, x) = \prod_{(U, \mathcal{Y}) \in \mathcal{F}(W^1,W^2)} Y_{U}(v, x) = \sum_{n \in \frac{1}{T_3}\Z} v_n^{\mathcal{S}} x^{-n-1}.
\]

\begin{remark}

\begin{enumerate}[{(1)}]
\item  Since different elements of $\mathcal{F}(W^1,W^2)$ can share the same underlying module $U$, the direct product $\mathcal{S}$ may contain many isomorphic copies of $U$.

\item By construction, each component module $U$ in $\mathcal{S}$ is a generalized $g_3$-twisted $V$-module, and hence the vertex operator $Y_{\mathcal{S}}(-, x)$ satisfies the twisted Borcherds identity (See \ref{tensor-product}) when applied componentwise. Consequently, $Y_{\mathcal{S}}(-, x)$ formally satisfies the twisted Borcherds identity on $\mathcal{S}$. However, unlike the direct sum case, the action on a direct product does not guarantee the lower-truncation property.
This absence causes the sums in the twisted Borcherds identity to become infinite. 
Therefore, $\left(\mathcal{S}, Y_{\mathcal{S}}(-, x)\right)$ is generally not a weak $g_3$-twisted $V$-module.

\item Let $W$ be a subspace of $\mathcal{S}$ that is invariant under the operators $v^{\mathcal{S}}_n$ for all $v \in V$.
If the restriction of $Y_{\mathcal{S}}(v , x)$ to $W$  is lower-truncated, then $W$ carries the structure of a weak $g_3$-twiseted $V$-module.

\end{enumerate}
\end{remark}

For $w^1 \in W^1$, $w^2 \in W^2$, $n \in \mathbb{C}$ and $k \in \mathbb{N}$, we define
\[
(w^1)^{\diamond}_{n;k}w^2 = \prod_{(U, \mathcal{Y}) \in \mathcal{F}(W^1,W^2)} \big( (w^1)^{\mathcal{Y}}_{n;k}w^2 \big) \in \mathcal{S}.
\]
Specifically, for each $(U,\mathcal{Y}) \in \mathcal{F}(W^1, W^2)$, the $(U,\mathcal{Y})$-th component of $(w^1)^{\diamond}_{n;k}w^2$ is exactly $(w^1)^{\mathcal{Y}}_{n;k}w^2$.
Furthermore, we define the formal series
\[
\mathcal{Y}^{\diamond}(w^1,x)w^2 = \sum_{\substack{n \in \mathbb{C} \\ k \in \mathbb{N}}} (w^1)^{\diamond}_{n;k}w^2 \, x^{-n-1} (\log x)^k.
\]

Let 
\[
W^1 \Diamond W^2=\text{span} \big\{ (w^1)^{\diamond}_{n;k}w^2 \  |  \  w^1 \in W^1, w^2 \in W^2, n \in \C, k \in \N \big\} \subset \mathcal{S}.
\]
Even though it has not yet been proven that $W^1 \Diamond W^2$ is a weak $g_3$-twisted $V$-module, 
we extend the terminology of Definition \ref{weight} to it, 
we will use the notations $(W^1 \Diamond W^2)_{(\lambda)}$, $\wt (w)$, and $\wt (W^1 \Diamond W^2)$ as in Definition \ref{weight}.

\vspace{0.5cm}
 
We collect the following basic facts.
\begin{enumerate}[{(1)}]
\item  Lemma \ref{U-Y-Lemma} implies that there exists a $K\in \N$ such that
 \[
 \mathcal{Y}^{\diamond}(w^1 ,x)w^2 \in \mathcal{S} \{x\} [\text{log}(x)]_{\leq K}.
 \]

\item By Lemmas \ref{U-Y-Lemma} and \ref{lem:s-n}, we have
   \[
	\big(L(0)_n\big)^d U = 0 \quad \text{for any } (U, \mathcal{Y}) \in \mathcal{F}(W^1,W^2).
	\]
In other words, for any such $(U, \mathcal{Y})$,
\[
\big( L(0) -(\wt w^1 + \wt w^2 -n-1) \big)^d \big( (w^1)^{\Y}_{n;k}w^2 \big)=0.
\]
Passing to the direct product, we obtain
\begin{align} \label{wt-6-1}
\big( L(0) -(\wt w^1 + \wt w^2 -n-1) \big)^d \big( (w^1)^{\diamond}_{n;k}w^2 \big)=0.
\end{align}
Hence, $W^1 \Diamond W^2$ decomposes into a direct sum of generalized eigenspaces of $L(0)$.
Moreover, by Equation \eqref{wt-6-1} and Lemma \ref{U-Y-Lemma}, it follows that 
\[
\wt (W^1 \Diamond W^2) \subset \bigcup^n_{i=1} \left(\lambda_i + \frac{1}{T_3} \mathbb{N} \right).
\]
We therefore conclude that the vertex operator $\Y^{\diamond}(-, x)$ is lower-truncated: 
for any $w^1 \in W^1$, $w^2 \in W^2$, and $h \in \C$,
$$(w^1)^{\diamond}_{h+n;k}w^2=0  \ \ \text{for} \ n \in \N \ \ \text{sufficiently large, independently of} \ k. $$

\item  Let $v \in V,  w^1 \in W^1$, $w^2 \in W^2$ be homogeneous elements, and let $n \in \C, p \in \frac{1}{T_3}\Z$, and $k \in \N$.
For any $(U, \mathcal{Y}) \in \mathcal{F}(W^1,W^2)$, the following weight identity holds:
\[
\wt \left(v_p\left((w^1)^\Y_{n; k}w^2 \right) \right) =\wt v -p-1 + \wt w^1 + \wt w^2 -n-1.
\]
This implies the same relation for the direct product operator:
\[
 \wt \big(  v^{\mathcal{S}}_p \left(  (w^1)^{\diamond}_{n;k}w^2 \right) \big)=\wt v -p-1 + \wt w^1 + \wt w^2 -n-1.
\]
Consequently, for $w \in W^1 \Diamond  W^2$, we have $v^{\mathcal{S}}_p w =0$ for $p \in \frac{1}{T_3}\Z$ sufficiently large.

\end{enumerate}
 
\begin{lemma}
$(W^1 \Diamond W^2,  Y_{\mathcal{S}}(-, x) )$ is a generalized $g_3$-twisted $V$-module, and
$\Y^{\diamond}$ is a twisted logarithmic intertwining operator of type $\binom{W^1 \Diamond W^2}{W^1 \ W^2}$.

\end{lemma}

\begin{proof}

Since the twisted operator $\Y^{\diamond}(-, x)$ is lower-truncated, the twisted Borcherds identity for its individual components implies that the identity holds for $\Y^{\diamond}(-, x)$ itself. By repeating the argument in Lemma \ref{ImY-1}, we find that $W^1 \Diamond W^2$ is closed under the action of $v^{\mathcal{S}}_p$ for all $v \in V$ and $p \in \frac{1}{T_3}\mathbb{Z}$. Thus, $W^1 \Diamond W^2$ forms a weak $g_3$-twisted $V$-module and is, consequently, a generalized $g_3$-twisted $V$-module.

Furthermore, as each component of $\Y^{\diamond}(-, x)$  satisfies the $L(-1)$-derivative property, so does  $\Y^{\diamond}(-, x)$ itself.
Consequently, $\Y^{\diamond}(-, x)$ is a twisted logarithmic intertwining operator.

The proof is complete.
\end{proof}

Furthermore, we have
\begin{theorem}
$( W^1 \Diamond W^2, \mathcal{Y}^{\diamond})$ satisfies the universal property of the fusion product of $W^1$ and $W^2$. Therefore,
$W^1 \Diamond W^2 \cong W^1 \boxtimes W^2.$
\end{theorem}

\begin{proof}
Let $(W, \mathcal{J})$ be an arbitrary twisted logarithmic intertwining operator of type $\binom{W}{W^1 W^2}$,
where $W$ is a generalized $g_3$-twisted $V$-module.
By the definition of $\mathcal{F}_{(W^1,W^2)}$, there exists a pair $(U, \mathcal{Y}) \in \mathcal{F}_{(W^1,W^2)}$ 
and a $g_3$-twisted $V$-module isomorphism $f: U \to \operatorname{Im} \mathcal{J}$ such that
$f \circ \mathcal{Y} = \mathcal{J}$.

Let
\[
\pi_{(U, \Y)}: \mathcal{S} \to U
\] 
be the canonical projection onto the $(U,\mathcal{Y})$-component, 
and let
\[
P_{(U, \Y)} = \pi_{(U, \Y)}|_{W^1 \Diamond W^2}: W^1 \Diamond W^2 \to U
\]
denote its restriction to the subspace $W^1 \Diamond W^2 \subseteq \mathcal{S}$. 

Define $h = f \circ P_{(U, \Y)}$. Then $h$ is a $g$-twisted $V$-module homomorphism and satisfies
\[
h \circ \mathcal{Y}^{\diamond} = \mathcal{J}.
\]
Since $\mathcal{Y}^{\diamond}$ is surjective,
the homomorphism $h$ satisfying $h \circ \mathcal{Y}^{\diamond} = \mathcal{J}$ is unique.

This establishes that $(W^1 \Diamond W^2, \mathcal{Y}^{\diamond})$ satisfies the universal property 
of the fusion product of $W^1$ and $W^2$. The proof is complete.

\end{proof}

\subsection{A Special Case} \label{sec6-2}

{\bf Throughout this subsection, we always assume that, up to isomorphism, $V$ has only finitely many irreducible grading restricted generalized $g_3$-twisted modules.}
Let $S^1, S^2, \cdots, S^p$ be a complete set of representatives of the isomorphism classes of irreducible grading restricted generalized $g_3$-twisted $V$-modules,
with conformal weights $u_1, u_2, \cdots, u_p$, respectively.
Each module decomposes as 
\[ 
S^i= \bigoplus_{n \in \N}S^i_{(u_i+ n/T_3)}
\] 
For a grading restricted generalized $g_3$-twisted $V$-module $W$, we define the natural number $d_W$ by 
\[
d_W=\text{dim} W_{(u_1)} + \cdots +\text{dim} W_{(u_p)}. 
\]
We remark that if $W$ is a grading restricted generalized $g_3$-twisted $V$-module, then $W$ is zero if and only if $d_W=0$.

\vspace{0.5cm}

Lemma \ref{leq-d} implies the existence of complex numbers $\lambda_1, \dots, \lambda_n$ with the following property: for every $(W, \Y)$ in $\mathcal{F}{(W^1,W^2)}$ and any homogeneous complementary subspace $E_W$ of $C_1(W)$ in $W$, we have
\[
 E_W \subseteq \bigoplus_{i=1}^n W_{(\lambda_i)}  \ \  \ \text{and} \ \ \ \dim (E_W) \leq m.
\]
Furthermore, by Proposition \ref{c1-fg}, there exist natural numbers $d_{u_1}, d_{u_2}, \cdots, d_{u_p}$ such that for any $(W, \Y) \in \mathcal{F}{(W^1,W^2)}$,
\[
\dim (W_{(u_i)}) \leq d_{u_i} \ \ \ \text{for} \ \ i=1, 2, \cdots, p.
\]
Let $d_0 = \sum_{i=1}^p d_{u_i}$. It follows that $d_W \leq d_0$ for any $(W, \Y) \in \mathcal{F}{(W^1,W^2)}$.

Assuming the existence of the fusion product 
$W^1 \boxtimes W^2$, then for any $(U, \Y) \in \mathcal{F}_{(W^1,W^2)}$, $U$ is a homomorphic image of $W^1 \boxtimes W^2$, and hence 
\[
d_U \leq d_{W^1 \boxtimes W^2}  \leq d_0.
\]
For this reason, we let $(W^1 \odot W^2, \Y^\odot) \in \mathcal{F}_{(W^1,W^2)}$ be such that $d_{W^1 \odot W^2}$ is maximal. 
We will show that this pair $(W^1 \odot W^2, \Y^\odot)$ satisfies the universal property of the fusion product of $W^1$ and $W^2$. 

\begin{theorem} \label{tensor-product}
The pair $(W^1 \odot W^2, \Y^{\odot})$ satisfies the universal property of the fusion product of $W^1$ and $W^2$. Therefore, 
\[
(W^1 \odot W^2 , \Y^\odot)\cong (W^1 \boxtimes W^2, \Y^{\boxtimes} ).
\]
\end{theorem}

\begin{proof}
Let $(U, \mathcal{J})$ be an arbitrary twisted logarithmic intertwining operator of type $\binom{U}{W^1 W^2}$.
To relate $U$ to  $W^1 \odot W^2$, we consider the direct sum intertwining operator
\[
\Y^{\odot} \oplus \mathcal{J} \in  \mathcal{I}\binom{(W^1 \odot W^2) \oplus U}{W^1 \ W^2}.
\]
Let 
\[
p\colon (W^1 \odot W^2) \oplus U \to W^1 \odot W^2  \ \ \text{and} \ \ q\colon (W^1 \odot W^2) \oplus U  \to U
\]
denote the canonical projections. 
Let $G = \mathrm{Im}(\Y^{\odot} \oplus \mathcal{J})$ and $\mathcal{Y}^* = \Y^{\odot} \oplus \mathcal{J}$. 
Then $(G, \mathcal{Y}^*) \in \mathcal{F}_{(W^1,W^2)}$, 
and the restriction $p|_G\colon G \to W^1 \odot W^2$ is an epimorphism.
It follows that 
\[
\dim G_{(u_i)} \geq \dim (W^1 \odot W^2)_{(u_i)}
\]
for any $i=1,2, \cdots, p$, and hence $d_{G} \geq d_{W^1 \odot W^2}$.
By the maximality of $d_{W^1 \odot W^2}$, we obtain 
\[
\dim G_{(u_i)} = \dim (W^1 \odot W^2)_{(u_i)}
\]
for any $i=1,2, \cdots, p$. It follows that $d_{\text{Ker}\ p|_{G} }=0$. 
Hence $p|_{G}$ is an isomorphism.

Let $\phi=q\circ p|_{G}^{-1} $, then $h$ is a generalized $g_3$-twisted $V$-module homomorphism and 
satisfies
\begin{align} \label{y=j}
\phi \Y^\odot(w^1, x) w^2 =\mathcal{J}(w^1, x)w^2
\end{align}
for any $w^1 \in W^2$ and $w^2 \in W^2.$
Since $\mathcal{Y}^\odot$ is surjective,
the homomorphism $\phi$ satisfying \ref{y=j} is unique.

This establishes that $(W^1 \odot W^2, \mathcal{Y}^{\odot})$ satisfies the universal property of the fusion product of $W^1$ and $W^2$. The proof is complete.
\end{proof}

The following tautological construction is well known in the untwisted case:
\begin{proposition}
Let $V$ be a $C_2$-cofinite vertex operator such that the category of generalized $g_3$-twisted $V$-modules is semi-simple. 
Let $S^1, S^2, \cdots, S^p$ be a complete set of representatives of the isomorphism classes of irreducible generalized $g_3$-twisted $V$-modules,
For $k=1,2$, let $W^k$ be a $C_1$-cofinite generalized $g_k$-twisted $V$-module. 
Then we have
\[
W^1 \boxtimes W^2 = \bigoplus_{i=1}^p (S^i)^{\oplus N(W^1, W^2; S^i)}
\]
\end{proposition}

\begin{proof}
By Theorem~\ref{tensor-product}, the fusion product $W^1 \boxtimes W^2$ exists. 
Since $W^1 \boxtimes W^2$ is a semi-simple module, it is completely determined up to isomorphism by the multiplicities of the simple modules $S^i$ in it.

The universal property of the tensor product yields the following vector space isomorphism for any $g_3$-twisted $V$-module $U$:
\[
\operatorname{Hom}_{V}(W^1 \boxtimes W^2, U) \cong \mathcal{I} \binom{U}{W^1 \; W^2}.
\]
This implies that $\dim \text{Hom}_{V}(S^i, W^1 \boxtimes W^2) = N(W^1, W^2; S^i)$ for all $i$.
Furthermore, by Theorem~\ref{finite-dim}, the multiplicity spaces $N(W^1, W^2; S^i)$ are finite-dimensional for each $i$.
We thus obtain the decomposition:
\[
W^1 \boxtimes W^2 = \bigoplus_{i=1}^p (S^i)^{\oplus N(W^1, W^2; S^i)}.
\]

This completes the proof.
\end{proof}





\end{document}